# Undecidable propositions with Diophantine form arisen from every axiom and every theorem in Peano Arithmetic

## T. Mei

(Department of Journal, Central China Normal University, Wuhan, Hubei PRO, People's Republic of China

E-Mail:   meitao@mail.ccnu.edu.cn     meitaowh@public.wh.hb.cn )

**Abstract:** By constructing a provability predicate with Diophantine form in Peano Arithmetic (PA) and based on the abstract form of Gödel's Second Incompleteness Theorem, it is proved that, if PA is consistent, then for every axiom and every theorem of PA, we can construct a corresponding undecidable proposition with Diophantine form. Finally, we present a method that transforms seeking a proof of a mathematical (set theoretical, number theoretical, algebraic, geometrical, topological, etc) proposition into solving a Diophantine equation.

Since the Matijasevič-Robinson-Davis-Putnam (MRDP) theorem has been established in the year 1970, some undecidable propositions with Diophantine form in First-order Peano Arithmetic (PA) have been studied[1, 2]. For example, T. Franzén[3] explains how the MRDP theorem can be used to obtain a proof to Gödel's First Incompleteness Theorem. Another consequence of the MRDP theorem mentioned in Ref.[3], which depends on the fact that PA (in fact, EA (Elementary Arithmetic)[4]) proves the MRDP theorem, is that, given any $\Pi^0_1$- sentence $K$, we can construct a corresponding Diophantine equation $D(x_1, \cdots, x_n)=0$ such that it is provable in PA that $K$ is true if and only if the equation has no solution. According to this consequence, if a $\Pi^0_1$- sentence is undecidable, then we can construct an undecidable proposition with Diophantine form.

In this paper, by constructing a provability predicate with Diophantine form in PA and based on the abstract form of Gödel's Second Incompleteness Theorem, it is proved that, for every axiom and every theorem of PA, we can construct a corresponding Diophantine equation such that if PA is consistent, then the Diophantine equation has no solution, but this fact is not provable within PA. Hence, every axiom and every theorem of PA can draw forth an undecidable proposition.

The method obtaining the above result is quite simple and don't need more knowledge about PA. For example, we don't need the conclusion "PA proves the MRDP theorem". This method may also be used to those arithmetic theories that do not prove the totality of exponentiation.

**Lemma 1.** (1) Introducing a set $\mathrm{Prf}_F$ for a formula $F$ of PA:

$\mathrm{Prf}_F := \{a \in \mathrm{Prf}_F \mid a$ is a Gödel code of a proof of the formula $F$ of PA$\}$,

the set $\mathrm{Prf}_F$  is recursive.

(2) There is a Diophantine equation $D_{\mathrm{PE}}(\#F, x_0, \{x_w\})=0$, which is called *the proof equation of a formula $F$ of PA*, such that $a \in \mathrm{Prf}_F$ if and only if $D_{\mathrm{PE}}(\#F, x_0, \{x_w\})=0$ has solution $x_0=a$ and $\{x_w\}=\{x_{w0}\}$, where $\#F$ reads the Gödel code of the formula $F$, $\{x_w\}$ is the abbreviation for $w$ variables $x_1, x_2, \cdots, x_w$. For the proof equation we have

(2.1) $D_{\text{PE}}(\#F, x_0, \{x_w\})=0$ has solution if and only if PA ⊢ $F$;

(2.2) If PA ⊢ $F$, then if PA is consistent, then $D_{\text{PE}}(\#(\neg F), x_0, \{x_w\})=0$ has not solution, where $\#(\neg F)$ reads the Gödel code of the formula $\neg F$.

(3) There is a formula $d_{\text{PE}}((\#F)', \boldsymbol{x_0}, \{x_w\})$ of PA as the representation of $D_{\text{PE}}(\#F, x_0, \{x_w\})=0$ such that, if $D_{\text{PE}}(\#F, b_0, \{b_w\})=0$, then PA ⊢ $d_{\text{PE}}((\#F)', \boldsymbol{b_0'}, \{\boldsymbol{b_w'}\})$; If $D_{\text{PE}}(\#F, b_0, \{b_w\})\neq0$, then PA ⊢ $\neg d_{\text{PE}}((\#F)', \boldsymbol{b_0'}, \{\boldsymbol{b_w'}\})$. In the formula $d_{\text{PE}}((\#F)', \boldsymbol{x_0}, \{x_w\})$, $(\#F)'$ is the abbreviation for term $\boldsymbol{0''\cdots'}$ (There are $\#F$ successors "'"), $\boldsymbol{b_0'}$ and $\{\boldsymbol{b_w'}\}$ as well.

(4) Let $B(\boldsymbol{x})$ be the formula $\exists(x_0, \{x_w\})\, d_{\text{PE}}(\boldsymbol{x}, x_0, \{x_w\})$. For which we have

(4.1) If $D_{\text{PE}}(\#F, x_0, \{x_w\})=0$ has solution, then PA ⊢ $B((\#F)')$;

(4.2) If PA ⊢ $B((\#F)')$, then if PA is $\omega$-consistent, then $D_{\text{PE}}(\#F, x_0, \{x_w\})=0$ has solution.

(5) Any formula $F$ of PA can cause a sequence of proof equations and $B$-formulas. Concretely, from the formula $F$ and the corresponding $D_{\text{PE}}(\#F, x_0, \{x_w\})=0$ and $B((\#F)')$, let $W_1$ be the $B((\#F)')$; For the formula $W_1$ of PA we have the corresponding $D_{\text{PE}}(\#W_1, x_{10}, \{x_{w1}\})=0$ and $B((\#W_1)')$, let $W_2$ be the $B((\#W_1)')$; For the formula $W_2$ of PA we have the corresponding $D_{\text{PE}}(\#W_2, x_{20}, \{x_{w2}\})=0$ and $B((\#W_2)')$; $\cdots$. The basic characteristic of such sequence of

$$D_{\text{PE}}(\#W_t, x_{t0}, \{x_{wt}\})=0 \text{ and } B((\#W_t)') \qquad (t=0, 1, 2, \cdots; W_{t+1}:= B((\#W_t)'), W_0:=F)$$

is that if PA is $\omega$-consistent, then if one equation $D_{\text{PE}}(\#W_u, x_{u0}, \{x_{wu}\})=0$ in the sequence has solution, then all the rest equations $D_{\text{PE}}(\#W_t, x_{t0}, \{x_{wt}\})=0$ have solutions and PA ⊢ $B((\#W_t)')$ for all $t=0, 1, 2, \cdots$.

The proof of Lemma 1 is quite simple, here we only write out that of Lemma 1(4.1) and (4.2) as two examples.

*Proof.* (4.1): Assuming a group of solution of $D_{\text{PE}}(\#F, x_0, \{x_w\})=0$ are $x_0=a$ and $\{x_w\}=\{a_w\}$, from Lemma 1(3) we have PA ⊢ $d_{\text{PE}}((\#F)', \boldsymbol{a'}, \{\boldsymbol{a_w'}\})$, and then PA ⊢ $\exists(x_0, \{x_w\})\, d_{\text{PE}}((\#F)', \boldsymbol{x_0}, \{x_w\})$, namely, PA ⊢ $B((\#F)')$.

(4.2): If PA ⊢ $B((\#F)')$, namely, PA ⊢ $\neg \forall(x_0, \{x_w\})\neg d_{\text{PE}}((\#F)', \boldsymbol{x_0}, \{x_w\})$. Suppose $D_{\text{PE}}(\#F, x_0, \{x_w\})=0$ has no solution, from Lemma 1(3), for all $(\boldsymbol{b_0'}, \{\boldsymbol{b_w'}\})$ we have PA ⊢ $\neg d_{\text{PE}}((\#F)', \boldsymbol{b_0'}, \{\boldsymbol{b_w'}\})$. Hence, all formulas $\neg d_{\text{PE}}((\#F)', \boldsymbol{b_0'}, \{\boldsymbol{b_w'}\})$ for all $(\boldsymbol{b_0'}, \{\boldsymbol{b_w'}\})$ and $\neg \forall(x_0, \{x_w\})\neg d_{\text{PE}}((\#F)', \boldsymbol{x_0}, \{x_w\})$ are provable in PA, PA is thus $\omega$-inconsistency.

**Remark 1.** (1) If we apply the MRDP theorem to the set $\text{Prf}_F$ introduced in Lemma 1(1), then what we obtain directly is a Diophantine equation $D_{\text{PE}}(\#F, a, \{x_w\})=0$ with the parameter $a$ such that $a \in \text{Prf}_F$ if and only if $D_{\text{PE}}(\#F, a, \{x_w\})=0$ has solution $\{x_w\}=\{x_{w0}\}$. This characteristic of the equation $D_{\text{PE}}(\#F, a, \{x_w\})=0$ allows us to replace $a$ with an unknown $x_0$ and, thus, obtain the corresponding proof equation $D_{\text{PE}}(\#F, x_0, \{x_w\})=0$.

(2) Changing the parameter $a$ in the equation $D_{\text{PE}}(\#F, a, \{x_w\})=0$ into an unknown $x_0$ leads to that great change take place in characteristics of the equation. Whether the equation $D_{\text{PE}}(\#F, a, \{x_w\})=0$ has solution for a given $a$ is decidable, since whether a given natural number is a Gödel code of a proof of a given formula of PA is decidable. Once we know a given $a$ is (not) a Gödel code of a proof of a given formula $F$, at the same time we also know that $D_{\text{PE}}(\#F, a, \{x_w\})=0$ has (not) solution. However, for the equation $D_{\text{PE}}(\#F, x_0, \{x_w\})=0$ we have not similar deciding method.

Of course, maybe we have different method to decide whether the equation $D_{\text{PE}}(\#F, x_0, \{x_w\})=0$ has solution. For example, if we know that $F$ is (not) a theorem of PA, then at the same time we also know that $D_{\mathrm{PE}}(\#F, x_0, \{x_w\})=0$ has (not) solution.

(3) If we rewrite the equation $D_{\mathrm{PE}}(\#F, x_0, \{x_w\})=0$ to the form $D_{\mathrm{PEL}}(\#F, x_0, \{x_w\})=D_{\mathrm{PER}}(\#F, x_0, \{x_w\})$, where both $D_{\mathrm{PEL}}(\#F, x_0, \{x_w\})$ and $D_{\mathrm{PER}}(\#F, x_0, \{x_w\})$ are polynomials with natural number coefficients, then from the form $D_{\mathrm{PEL}}(\#F, x_0, \{x_w\})=D_{\mathrm{PER}}(\#F, x_0, \{x_w\})$ we can obtain the corresponding representation $d_{\mathrm{PE}}((\#F)', \boldsymbol{x}_0, \{x_w\})$ immediately.

For obtaining $d_{\mathrm{PEL}}((\#F)', \boldsymbol{x}_0, \{x_w\})$ from $D_{\mathrm{PEL}}(\#F, x_0, \{x_w\})=D_{\mathrm{PER}}(\#F, x_0, \{x_w\})$, what we need is only the representability of recursive functions (in fact, that of polynomial functions), but not the conclusion "PA proves the MRDP theorem".

(4) In the proof of Lemma 1(4.2), what we need is in fact the $n$-fold form of $\omega$-consistency, so called "$\omega^n$-consistency[2]", since $d_{\mathrm{PE}}((\#F)', \boldsymbol{x}_0, \{x_w\})$ is multivariate but not one free variable. However, according to the lemma in Ref.[2, §2]: "*Suppose $R_0 \subseteq T$. If $T$ is $\omega$-consistency, then $T$ is $\omega^k$-consistency, for all $k$*", we can only use the condition of $\omega$-consistency in the proof of Lemma 1(4.2).

(5) It is obvious that the condition of "$\omega$-consistency" in Lemma 1(4.2) can be weakened to that of "1-consistency[5]". However, we don't discuss "1-consistency" in this paper.

(6) All the conclusions of Lemma 1(3), (4.1) and (4.2) hold not only for $D_{\mathrm{PE}}(\#F, x_0, \{x_w\})=0$ and $d_{\mathrm{PE}}((\#F)', \boldsymbol{x}_0, \{x_w\})$, but also for any Diophantine equation $D(c, \{x_m\})=0$ and the corresponding representation $d_{\mathrm{D}}(\boldsymbol{c}', \{x_m\})$ in PA. ∎

Given different formulas of PA, maybe the forms of the corresponding proof equations and the representations in PA are different. For instance, maybe the forms of $D_{\mathrm{PE}}(\#F, x_0, \{x_u\})=0$ is different from that of $D_{\mathrm{PE}}(\#(\neg\, F), x_0, \{x_v\})=0$, the corresponding representations $d_{\mathrm{PE}}((\#F)', \boldsymbol{x}_0, \{x_w\})$ and $d_{\mathrm{PE}}((\#(\neg\, F))', \boldsymbol{x}_0, \{x_w\})$ as well. On the other hand, we have the following result.

**Lemma 2.** (1) There exists *an universal proof equation* $D_{\mathrm{PE\text{-}PA}}(k, x_0, \{x_n\})=0$ *for PA* such that once we replace $k$ with the Gödel code $\#F$ of a formula $F$ of PA, $D_{\mathrm{PE\text{-}PA}}(\#F, x_0, \{x_n\})=0$ becomes a proof equation of the formula $F$.

(2) The Diophantine equation $D_{\mathrm{PE\text{-}PA}}(k, x_0, \{x_w\})=0$ is undecidable.

**Remark 2.** (1) The pure formal proof of Lemma 2(1) is tedious; however, the basic idea of the proof is actually quite simple. We explain it as follows.

The rule of Gödelian coding insures that whether a given number $a$ is a Gödel code of a proof of a given formula $F$ of PA is decidable, we can design a computer program to do this thing. Such program is demanded to perform three tasks:

Task 1: Checking whether $a$ expresses a sequence of formulas $[F_1, F_2, \cdots, F_n]$ according to Gödelian coding rule.

Task 2: Checking whether a sequence of formulas $[F_1, F_2, \cdots, F_n]$ is a proof of the last formula $F_n$ in the sequence. Concretely, the program checks whether every formula $F_i (1 \leqslant i \leqslant n)$ ① is an axiom, or ② follows from $F_l$ and $F_m$ ($l < i, m < i$) by MP rule, or ③ follows from $F_j$ ($j < i$) by GEN rule.

By the way, we can improve slightly the rule of Gödelian coding of a proof of a formula by prescribing that, in the Gödel code of a proof of a formula $F$, every formula $F_i (1 \leqslant i \leqslant n)$ must take with its "warrantor". Concretely, ① If $F_i$ is an axiom, then it brings with the corresponding "axiom schema number" (We can assign a number for every axiom schema in advance); ② If $F_i$ follows from $F_l$ and $F_m$ ($l < i, m < i$) by MP rule, then it brings with the numbers $l$ and $m$; ③ If $F_i$ follows from $F_j$ ($j < i$) by GEN rule, then it brings with the number $j$. Designing such a coding system is very easy, but it brings great convenience to the design of the program.

Task 3: Checking whether the last formula $F_n$ in a sequence of formulas $[F_1, F_2, \cdots, F_n]$ is the given formula $F$.

Notice that both Task 1 and Task 2 are independent of the given formula $F$, hence, if the process of the program execution is Task 1 $\xrightarrow{\text{YES}}$ Task 2 $\xrightarrow{\text{YES}}$ Task 3, then this program is available for any formula in PA.

We therefore obtain an universal program, or an universal algorithm, $C(k, a)$. Once $k$ in $C(k, a)$ is replaced by the Gödel code $\#F$ of a formula $F$ in PA, $C(\#F, a)$ becomes a program checking whether a number $a$ is a Gödel code of a proof of the formula $F$.

Applying the MRDP theorem to the algorithm $C(k, a)$, we can find out a corresponding Diophantine equation $D_{\text{PE-PA}}(k, a, \{x_w\})=0$ with two parameters $k$ and $a$ such that once $k$ is replaced by the Gödel code $\#F$ of a formula $F$ in PA, $D_{\text{PE-PA}}(\#F, a, \{x_w\})=0$ has solution if and only if $a$ is a Gödel code of a proof of the formula $F$.

And then, by replacing $a$ in $D_{\text{PE-PA}}(k, a, \{x_w\})=0$ with an unknown $x_0$, we obtain an universal proof equation $D_{\text{PE-PA}}(k, x_0, \{x_n\})=0$ for any formula of PA.

(2) It is obvious that the method obtaining an universal proof equation for PA presented in Remark 2(1) is constructive. Of course, how to construct a *simple* universal proof equation for PA is a different matter.

(3) As well-known, whether a formula of PA is a theorem in PA is undecidable. Hence, for the set FC: $= \{b \in \text{FC} \mid b$ is a Gödel code of a formula of PA$\}$, if $k \in \text{FC}$, then whether the equation $D_{\text{PE-PA}}(k, x_0, \{x_w\})=0$ has solution is undecidable. (Notice that by an analysis similar to Remark 1(2) we see that if $k \in \text{FC}$ then whether the equation $D_{\text{PE-PA}}(k, a, \{x_w\})=0$ has solution for a given $a$ is decidable.) Now that whether the equation $D_{\text{PE-PA}}(k, x_0, \{x_w\})=0$ has solution is undecidable for $k \in \text{FC}$, which is a subset of the set of natural numbers, $D_{\text{PE-PA}}(k, x_0, \{x_w\})=0$ is thus undecidable and can be as an example of undecidable Diophantine equations. ∎

Once we obtain an universal proof equation $D_{\text{PE-PA}}(k, x_0, \{x_n\})=0$ for PA, for which we can introduce the corresponding representation $d_{\text{PE-PA}}(\boldsymbol{k}, \boldsymbol{x}_0, \{\boldsymbol{x}_n\})$ in PA, and by $B_{\text{PA}}(\boldsymbol{x})$ we denote the formula $\exists(\boldsymbol{x}_0, \{\boldsymbol{x}_w\})\, d_{\text{PE-PA}}(\boldsymbol{x}, \boldsymbol{x}_0, \{\boldsymbol{x}_w\})$.

All the conclusions in Lemma 1 hold for $D_{\text{PE-PA}}(k, x_0, \{x_n\})=0$, $d_{\text{PE-PA}}(\boldsymbol{k}, \boldsymbol{x}_0, \{\boldsymbol{x}_n\})$ and $B_{\text{PA}}(\boldsymbol{x})$. Especially, immediately from Lemma 1(2.1) and Lemma 1(4.1) we obtain

(P1)　If PA $\vdash A$, then PA $\vdash B_{\text{PA}}((\#A)')$.

And, further, we assume that $B_{\text{PA}}(\boldsymbol{x})$ satisfies the following two properties.

**Hypothesis**. For arbitrary formulas $A$, $A_1$ and $A_2$ of PA,

(P2)　PA $\vdash B_{\text{PA}}((\#(A_1 \rightarrow A_2))') \rightarrow (B_{\text{PA}}((\#A_1)') \rightarrow B_{\text{PA}}((\#A_2)'))$.

(P3)　PA $\vdash B_{\text{PA}}((\#A)') \rightarrow B_{\text{PA}}((\#(B_{\text{PA}}((\#A)')))')$.

As well-known, (P1) ∼ (P3) are so called Hilbert-Bernays-Löb derivability conditions[6, 7, 8], we are not going to discuss the proof of (P2) and (P3) here and, thus, regard them as "Hypothesis".

Generally speaking, if PA is $\omega$-consistent, then $B_{\text{PA}}(\boldsymbol{x})$ is stronger than the provability predicate $Proof_F$ in a system $T$. For example, apart from satisfying (P1) ∼ (P3), for $B_{\text{PA}}(\boldsymbol{x})$ we have "If PA $\vdash B_{\text{PA}}((\#F)')$, then PA $\vdash F$" (Immediately from Lemma 1(4.2) and Lemma 1(2.1)). However, the provability predicate $Proof_F$ in $T$ has not such property when $T$ is not demanded to be $\omega$-consistent[6]. And, in the discussion below, we also do not demand that PA is $\omega$-consistent.

If (P1) ∼ (P3) hold for $B_{\text{PA}}(\boldsymbol{x})$, then we have

**Lemma 3 (Löb's Theorem)**. For any formula $A$, if PA $\vdash B_{\text{PA}}((\#A)') \rightarrow A$, then PA $\vdash A$.

**Lemma 4 (Gödel's Second Incompleteness Theorem, abstract form).** If PA is consistent, then if $A$ is an axiom or a theorem of PA, then not PA $\vdash \neg\, B_{PA}((\#(\neg\, A))')$.

**Remark 3.** The difference between Lemma 3 and the standard form of the Löb's Theorem (See, for example, Refs. [6] or [8]) is only that the provability predicate *Proof$_F$* in the standard form is replaced by $B_{PA}(x)$ given by this paper, the difference between Lemma 4 and the standard form of the abstract form of Gödel's Second Incompleteness Theorem as well. We therefore can follow those proofs of the standard forms in literatures to prove Lemma 3 and Lemma 4. For example, we can copy verbatim the proof of Theorem 18.4 in Ref. [6] as that of Lemma 3, as long as in which $B(x)$ is replaced by the formula $B_{PA}(x)$. Similarly, we can obtain a proof of Lemma 4 through replacing the formula $0=1$ with $\neg\, A$, where $A$ is any axiom or theorem of PA, and replacing $B(x)$ with the formula $B_{PA}(x)$ in the proof of Theorem 18.3, in which the Löb's Theorem is used, in Ref. [6]. ∎

We now investigate what result can be obtained from Lemmas 1 ~ 4. If $A$ is an axiom or a theorem of PA, then according to Lemma 1(2.2), if PA is consistent, then the assertion "the Diophantine equation $D_{PE\text{-}PA}(\#(\neg\, A), x_0, \{x_w\})=0$ has no solution" is true. On the other hand, the corresponding representation in PA of this assertion is just $\neg\, B_{PA}((\#(\neg\, A))')$, and Lemma 4 shows that if PA is consistent, then $\neg\, B_{PA}((\#(\neg\, A))')$ is not provable in PA. We therefore obtain the following conclusion:

**Theorem.** For every axiom and every theorem of PA, we can construct a corresponding Diophantine equation such that if PA is consistent, then the Diophantine equation has no solution, but this fact is not provable within PA.

The above Theorem shows that undecidable propositions are at least as much as the summation of whole axioms and theorems in PA, since for every axiom and every theorem of PA, we can construct a corresponding undecidable proposition.

Using the proof equation of a formula of PA and the corresponding representation, we can do other things. For example, applying the Diagonal Theorem[6, 7, 8] to two formulas $\neg\, B_{PA}(x)$ and $B_{PA}(x)$, we obtain the Gödel sentence $G \leftrightarrow \neg\, B_{PA}((\#G)')$ and the Henkin sentence $H \leftrightarrow B_{PA}((\#H)')$ with Diophantine form, respectively.

On the other hand, given a proposition $f$ (or its negation) in elementary number theory, e.g., twin prime conjecture, Goldbach's conjecture, the Riemann hypothesis[9], P$\neq$NP[10], we can find out the corresponding formula $F$ (or $\neg\, F$) of PA. And then, we construct the corresponding proof equation $D_{PE}(\#F, x_0, \{x_w\})=0$ or $D_{PE\text{-}PA}(\#F, x_0, \{x_n\})=0$ of the formula $F$. Once we obtain a solution of $D_{PE}(\#F, x_0, \{x_w\})=0$ or $D_{PE\text{-}PA}(\#F, x_0, \{x_n\})=0$, at the same time we obtain a proof of the proposition $f$. And, if we prove that the equation $D_{PE}(\#F, x_0, \{x_w\})=0$ or $D_{PE\text{-}PA}(\#F, x_0, \{x_n\})=0$ has solution, then we know that the proposition $f$ holds. We can deal with the negation of $f$ in the same way.

We can generate the above method obtaining a proof of a proposition in elementary number theory to any mathematical (set theoretical, number theoretical, algebraic, geometrical, topological, etc) proposition. Namely, taking advantage of the idea of proof equation, maybe we can transform seeking a proof of a mathematical proposition into solving a Diophantine equation.

For realizing such transformation, given a proposition $f$ in a mathematical system $s$, we first construct a corresponding formal system $S$. The construction of $S$ can be quite "random". For example, we can add some known theorems as axioms, and try to use second-order logic, etc, because our purpose is now not to study a formal system. However, we make demands on $S$: ①

There is a formula $F$ of $S$ to express the proposition $f$ of $s$. ② Gödelian coding must be implemented in $S$.

And then, we introduce a set $\mathrm{Prf}_F$ for the formula $F$ of $S$:

$$\mathrm{Prf}_F := \{a \in \mathrm{Prf}_F \mid a \text{ is a Gödel code of a proof of the formula } F \text{ of } S\},$$

by applying the MRDP theorem to the recursive set $\mathrm{Prf}_F$ we find out a Diophantine equation $D_{\mathrm{PE}}(\#F, a, \{x_w\})=0$ such that $a \in \mathrm{Prf}_F$ if and only if $D_{\mathrm{PE}}(\#F, a, \{x_w\})=0$ has solution.

Finally, we replace $a$ in the equation $D_{\mathrm{PE}}(\#F, a, \{x_w\})=0$ with an unknown $x_0$ and try to solve the equation $D_{\mathrm{PE}}(\#F, x_0, \{x_w\})=0$. If we obtain a solution, then we obtain a proof of the original proposition $f$.

We emphasize again the differences between two equations $D_{\mathrm{PE}}(\#F, a, \{x_w\})=0$ and $D_{\mathrm{PE}}(\#F, x_0, \{x_w\})=0$ for formula $F$ of formal system $S$ (see Remark 1(2)). Hence, even if we have had a corresponding equation $D_{\mathrm{PE}}(\#F, x_0, \{x_w\})=0$ of a proposition $f$, we still need creative intellectual work of human being.